# A Two-Tier Algebraic Schema to Map $(A^4 - C^4)/(D^4 - B^4)$ onto the Natural Numbers

Paul A. Roediger (proediger@comcast.net)

**Abstract**. A brief history and two formulations of the Diophantine problem's requirements are presented. One tier consisting of three two-parameter solutions is studied for its ability to provide examples for the small natural numbers $a \in N$ considered. Nested within it is a second tier consisting of five shifted-square solutions of the form $a = u^2 + c$, where $u, c \in Q$. All told, they provide numerical examples for all but two $a \in N[1000]$, the set of natural numbers less than or equal to 1000. A few open questions remain. Does this scheme of solutions cover every $a \in N[1000]$? If so, might they account for all $a \in N$? Are the three $tier_1$ solutions redundant with respect to the $a$'s they provide? Do other $tier_1$ and shifted-square $tier_2$ solutions exist?

**Brief History.** Several authors report that numerical solutions can be found for practically any natural number $a$ one considers. Choudhry (1995) tabulated solutions for seventy-five $a \in N[101]$, remarking that the equation seemed to have solutions for all $a \in N$. That table was completed by Tito Piezas in 2013, who conjectured that non-trivial solutions exist for every $a \in N$. Tomita (2016) extended the tables to $a \in N[1000]$ with one exception $a = 967$. Andrew Bremner, as reported by Piezas (2015) and Choudhry (2017), contributed the elusive example. The truth of the conjecture was further reinforced in Tomita (2017) who tabulated solutions for all but two $a \in N[20000]$. Noam Elkies subsequently contributed examples for the two missing $a$'s, 9719 and 16329.

Zajta (1983) showed there are several one-parameter solutions for $a = 1$. Choudhry (1998) provided a parametric solution for $a = 4$ and Roediger (1972, 2015), provided parametric solutions for $a = 4$ and $a = 9$.

Euler (1780), Hayashi (1911) and Piezas (2015) note a shifted-square solution $a = u^2 - 3$. It yields, for $u = 7/2^2$, a solution for $a = 1$. Grigorief, according to Dickson (1908, p 647), also notes an $a = 2$ solution for $u = 1871/33^2$. Zajta (1983) and Choudhry (2017)] report an $a = u^2 + 2$ solution, which, for $u = 7/2^2$ provides another $a = 1$ example. Five shifted-square solutions having the form $a = u^2 + c$ are identified in Roediger (1972, 2015). In addition to the $c = -3$ and $c = 2$ solutions just mentioned, these others are $c = -1$, $c = 3$ and $c = 9/4$.

This note presents two formulations that lead to a two-tiered scheme consisting of three two-parameter solutions which includes a second tier of special cases consisting of the five known shifted-square solutions. It also reports on the number of numerical examples they provide within the modest $a \in N[1000]$ domain.

**Two Formulations.** As usual, substitute $A = p+q$, $C = p-q$, $D = r+s$ and $B = r-s$ in $A^4 + aB^4 = C^4 + aD^4$, followed by $p = ry$, $s = qx$ and $r = qt$. This leads to the requirement

(A)
$$\frac{ax^3 - y}{y^3 - ax} = t^2$$

Following Roediger (1972) and (2015), a linearized version of (A) is obtained with $\frac{x - \rho y}{y - a\rho x} = t^2$, where $\rho = \frac{xy+1}{ax^2 + y^2}$. Upon elimination of $x$ and $y$ with $x = \frac{t^2 + \rho}{\omega}$, $y = \frac{a\rho t^2 + 1}{\omega}$, requirement A becomes

(B)
$$a^2 \rho^3 t^4 + (3a\rho^2 - 1) t^2 + a\rho^3 = \omega^2$$

A variant of (B) was obtained by a different method in Izadi and Baghalaghdam (2017). There, advanced elliptic curve theory was applied to develop numerical solutions for various $a \in N$, and the conjectured domain of the solution space was extended to all $a \in Q$. Inspection of ($A$) and ($B$) reveals that $(ak^{-4}, xk, yk^{-1}, tk^{-1})$ and $(ak^{-4}, \rho k^2, tk, \omega k)$ are solutions if $(a, x, y, t)$ and $(a, \rho, t, \omega)$ are solutions, respectively. Some choices of the constant $k$ can reduce the number of variables involved, e.g., setting $\rho = \square$ (say $k^{-2}$) removes $\rho$ from the mix.

**Two-Parameter Nests.** Three restrictions are considered in turn: $y = ax$ in ($A$); $\rho = \square$ in ($B$); and $\rho(\rho+4) = \square$ in ($B$). Each restriction leads to a rational two-parameter nest which is tabulated below in Table 1.

**Table 1**. Special cases of $A$ and $B$ provide the first tier of two-parameter nests for $a$

| $Tier_1$ | 1st Restriction | $a$ | $p$ | $q$ | $r$ | $s$ |
|---|---|---|---|---|---|---|
| $A_1$ | $y = ax$ <br> $x = (u+v)/(u-v)$ <br> $t = (uv-1)/(u-v)$ | $\dfrac{(u-v)(uv+1)}{(u+v)(uv-1)}$ | $uv+1$ | $u-v$ | $uv-1$ | $u+v$ |
| $B_1$ | $\rho = 1,\ t = v$ <br> $w = av^2 - u$ | $\dfrac{u^2 + v^2}{(2u+3)v^2 + 1}$ | $u(au^2 + 1)$ | $au^2 - v$ | $u(au^2 - v)$ | $u^2 + 1$ |
| $B_2$ | $t = 1,\ z = a\rho$ <br> $\rho = (u-1)^2 / u$ <br> $w = (u(u-1)z - 1)v$ | $\dfrac{u(u-1)(u^2 + v^2)}{(uv^2 - 1)}$ | $u(v^2 + 1)$ | $v(u^2 - 1)$ | $v(u+1)$ | $v^2 - 1$ |

**Note**. $B_1$ and $B_2$ are related: $B_2(u,v) = B_1(u + v\sqrt{k}, \sqrt{k})$, for $u, v, k \in Q$ and satisfies $ku^2 + (3k+1)u - k^2 = (kv)^2$

**Shifted-Square Nests.** One more restriction provides a second tier of "shifted-square" nests for $a$, providing a simpler search domain for some otherwise hard to reach $a \in N[1000]$ examples. They are tabulated in the following table.

**Table 2**. Shifted-square solutions provide the second tier of one-parameter nests for $a$

| $Tier_2$ | 2nd Restriction | $a$ | $p$ | $q$ | $r$ | $s$ |
|---|---|---|---|---|---|---|
| $A_{11}$ | $v = (u^2 - 2)/u$ | $u^2 - 3$ | $u(u^2 - 3)$ | $2(u^2 - 1)$ | $u(u^2 - 1)$ | $2$ |
| $B_{11}$ | $\alpha = \dfrac{-2}{u^2},\ t = \dfrac{1}{u}$ | $u^2 - 1$ | $u(u^4 - 1)$ | $2u^2 - 1$ | $u^4 - u^2 + 1$ | $u(2u^2 - 1)$ |
| $B_{12}$ | $\alpha = \dfrac{1}{u^2},\ t = \dfrac{1}{u}$ | $u^2 + 2$ | $u(u^2 + 2)$ | $1$ | $u^2 + 1$ | $u$ |
| $B_{13}$ | $\alpha = \dfrac{3u^2 + 4}{u^2(u^2 + 2)}$ <br> $t = \dfrac{u}{u^2 + 2}$ | $u^2 + 2$ | $u(u^2 + 1)$ | $3(u^2 + 2)$ | $u(u^2 + 4)$ | $3$ |
| $B_{14}$ | $\alpha = 0,\ t = \dfrac{1}{u}$ | $u^2 + 3$ | $u(u^2 + 1)(u^2 + 3)$ | $1$ | $u^4 + 3u^2 + 1$ | $u$ |
| $B_{15}$ | $\alpha = -\dfrac{3}{2},\ t = u$ | $u^2 + \dfrac{9}{4}$ | $u\left(u^4 + \dfrac{9}{4}u^2 + 1\right)$ | $u^4 + \dfrac{9}{4}u^2 + \dfrac{3}{2}$ | $u\left(u^4 + \dfrac{9}{4}u^2 + \dfrac{3}{2}\right)$ | $u^2 + 1$ |
| $B_{21}$ | $u = \dfrac{v^2}{v^2 - 1}$ | $v^2 - 1$ | $v(v^4 - 1)$ | $2v^2 - 1$ | $v^4 - v^2 + 1$ | $v(2v^2 - 1)$ |

**Note**. $B_{11} \equiv B_{21}$

Figure 1 depicts relationships among the two-tier scheme of nests described in Table 1 and 2. The universe, $a \in N[1000]$, is represented by the $A \times B$ oval. Computer searches performed for this study found numerical representations for all but two $a$'s belonging to the fallout subset $(A_1 \cup B_1 \cup B_2)^C$, which is designated by $F$.

A Venn diagram containing counts for the eight $tier_1$ subsets is depicted in Figure 1.

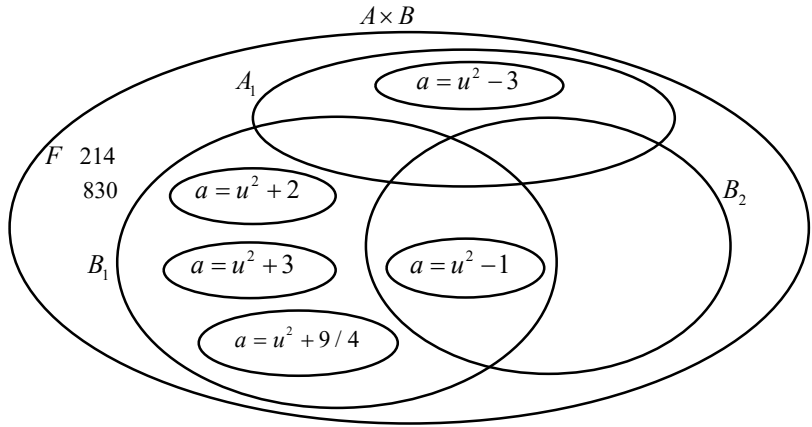

**Figure 1**. Nesting relationships among the various algebraic solution subsets defining an unsigned $a$

Search details and a coverage count summary are presented in Table 3 and summarized in Figure 2 below.

**Table 3**: Count summary of unique $a \in N[1000]$ provided by solutions $A_1$, $B_1$, $B_2$ and their various subsets

| Restriction | | $a$ | $M$ * | Direct ($D$) | Indirect ($I$) | $D+I$ | Tot | Cum |
|---|---|---|---|---|---|---|---|---|
| $A_1$ | | $\dfrac{(u-v)(uv+1)}{(u+v)(uv-1)}$ | 750 | 745 | 89 | 834 | | 834 |
| | $A_{11}$ | $u^2 - 3$ | 827 | 181 | 30 | 211 | 855 | 855 |
| $B_1$ | | $\dfrac{u^2 + v^2}{(2u+3)v^2 + 1}$ | 450 | 822 | 92 | 914 | | 981 |
| | $B_{11}$ | $u^2 - 1$ | 1383 | 733 | 59 | 791 | | 993 |
| | $B_{12}$, $B_{13}$ | $u^2 + 2$ | 707 | 155 | 26 | 181 | | 997 |
| | $B_{14}$ | $u^2 + 3$ | 80 | 96 | 13 | 109 | | 998 |
| | $B_{15}$ | $u^2 + \dfrac{9}{4}$ | 113 | 43 | 13 | 56 | 983 | 998 |
| $B_2$ | | $\dfrac{u(u-1)(u^2+v^2)}{(uv^2-1)}$ | 350 | 456 | 81 | 537 | | 998 |
| | $B_{21}$ | $v^2 - 1$ | 1383 | 733 | 59 | 791 | 862 | 998 |

* **Notes**: Search domains cover rational $u, v$ such that $1 < \max(h(u), h(v)) \leq M$, where $h$ is a height function of a fractional argument given by $h(i/j) = \max(|i|, |j|)$. Indirect $a$'s are direct $a$'s multiplied by $r^2$, all stemming from the fact that $A^4 + aB^4 = C^4 + aD^4$ is equivalent to $(rB)^4 + ar^2 A^4 = (rD)^4 + ar^2 C^4$, for all rational $r$.

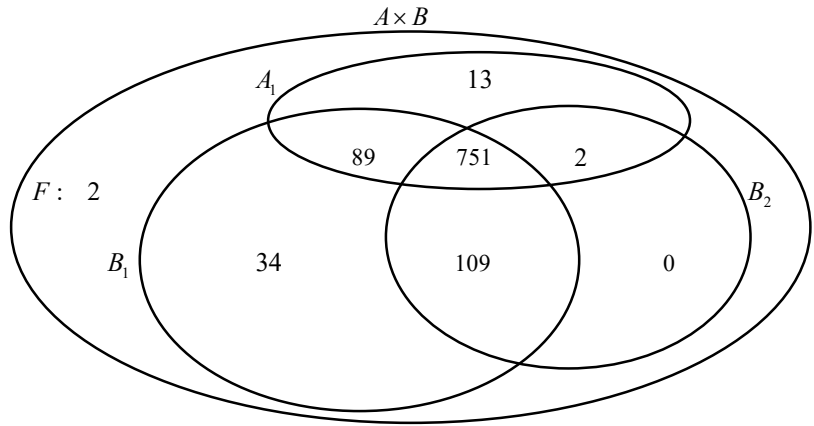

**Figure 2**. Counts of $a$'s $\in N[1000]$ covered by the seven $tier_1$ subsets plus the fallout subset $F$

**Search Wrap-up**. Having reached a practical hard stop for simply doing more $tier_1$ searching (and incurring prohibitively long run times with low expectation for getting any new $a$'s), attention was turned to further probe $tier_2$ solutions in the two subsets that stand out in Figure 3, i.e., the two $a$'s in $\{263, 670\}$, constituting the $A_1 \cap B_2 \cap B_1^C$ subset, and the two $a$'s in $\{214, 830\}$, constituting the fallout subset $F$.

Since all four $a$'s have a prime factor congruent to three modulo four, searching for a $u^2 + 9/4$ representation could be dispensed with immediately. Searches for representations among the other $u^2 + c$ formats were then simplified by examining the left-hand side (lhs) and right-hand side (rhs) of $u^2 + cX^4 = aY^4$ (modulo eight), for all combinations of $a \times c$ and all odd/even combinations of $u \times X \times Y$, where $a \in \{263, 670\}$ and $c \in \{-1, +2, +3\}$, and where $a \in \{214, 830\}$ and $c \in \{-3, -1, +2, +3\}$. Only those combinations whose $lhs$ and $rhs$ equated required the additional searches. The favorable searched cases are enumerated in the following table.

**Table 4**. Final searches exploring the possibility of giggling the counts presented in Figure 3

| a | c | u | X | Y | lhs | rhs | a | c | u | X | Y | lhs | rhs |
|---|---|---|---|---|---|---|---|---|---|---|---|---|---|
| 263 | −1 | Odd | Odd | Even | 0 | 0 | 214 | −3 | Odd | Odd | Odd | 6 | 6 |
|  | −1 | Even | Odd | Odd | 0 | 0 |  | −1 | Odd | Odd | Even | 0 | 0 |
|  | 3 | Even | Odd | Odd | 4 | 4 |  | 2 | Even | Odd | Odd | 3 | 3 |
| −263 | −1 | Odd | Odd | Even | 0 | 0 | −214 | −1 | Odd | Odd | Even | 0 | 0 |
|  | −1 | Odd | Even | Odd | 0 | 0 | 830 | −3 | Odd | Odd | Odd | 6 | 6 |
| 670 | −1 | Odd | Odd | Even | 0 | 0 |  | −1 | Odd | Odd | Even | 0 | 0 |
|  | 2 | Even | Odd | Odd | 3 | 3 |  | 2 | Even | Odd | Odd | 3 | 3 |
| −670 | −1 | Odd | Odd | Even | 0 | 0 | −830 | −1 | Odd | Odd | Even | 0 | 0 |

**Note**. All targeted searches were conducted over the entire range: $1 \le X \le 20000$, $1 \le Y \le rX$, where $r = (|a/c|)^{1/4}$

The searches described in Table 4 produced no new solutions. Two tentative conclusions may be drawn regarding these subsets: there are no $tier_2$ solutions moving 263 or 670 into the $B_1$ subset; nor are there any $tier_2$ solutions removing 214 or 830 from the fallout subset $F$.

**Conclusions**. In this computational study a credible case is made that the fallout subset $F$ is non-empty and possesses at most two members in $a \in N[1000]$.

One may ask:

(1) Are 214 and 830 representable in this scheme?
(2) Is it true that $B_2 \subset A \cup B_1$?
(3) Are there other $tier_1$ solutions?
(4) Are there other $tier_2$ shifted square solutions? (e.g., the discovery of an $a = u^2 - 9/4$ solution would cover the fallout subset $F$ since $a = 214 \times (2/5)^4$ and $a = 830 \times (3/11)^4$ for $u = 139/50$ and $u = 643/242$, respectively)
(5) Could such a schema be extended to include all natural numbers $a \in N$?

All searches and computations were performed on three laptops using R Statistical Software (v4.1.2 R CORE TEAM 2021), specifically *Rmpf*, an R-package by Maechler (2021). R programs (source code) enabling this work, along with lists of various numerical representations of $a$'s found, $a \in N[1000]$, are available upon request.

**References**. (All hyperlinks were successfully accessed on 3/28/2025)

Choudhry, A. (1995), "On the Diophantine equation $A^4 + hB^4 = C^4 + hD^4$", *Indian Journal of Pure and Applied Mathematics*, 26, 1057–1061

Choudhry, A. (1998), "The Diophantine Equation $A^4 + 4B^4 = C^4 + 4D^4$", *Indian Journal of Pure and Applied Mathematics*, 29 (11), 1127-1128, https://insa.nic.in/writereaddata/UpLoadedFiles/IJPAM/20005a05_1127.pdf

Choudhry, A. (2017), "A Note on the quartic Diophantine equation $A^4 + hB^4 = C^4 + hD^4$", *Notes on Number Theory and Discrete Mathematics*, 23(1), 1-3, https://nntdm.net/papers/nntdm-23/NNTDM-23-1-01-03.pdf

Dickson, L, (1920) History of the Theory of Numbers, Volume II
https://archive.org/details/historyoftheoryo02dickuoft/page/n3/mode/2up

Euler, L. (1780) Opera Arithmetica, Supplementum in Comment XXXIII T. I, 452-454,

Hayashi, T. (1911), "On the Diophantine Equation $X^4 + Y^4 = Z^4 + T^4$", *Tohoku Mathematics Journal*, 1, 143-145

Izadi, F. and Baghalaghdam, M. (2017), "Is the Quartic Diophantine Equation $A^4 + hB^4 = C^4 + hD^4$ Solvable for any Integer $h$?", *arXiv:1701.02602v3*

Mordell, L. J. (1969), Diophantine Equations, Academic Press, New York and London

Maechler, M. (2021), Rmpfr: R MPFR - Multiple Precision Floating-Point Reliable, R package version 0.8-4, https://CRAN.R-project.org/package=Rmpfr

Piezas, T. (2015), "On $a^4 + nb^4 = c^4 + nd^4$ and Chebyshev polynomials", available at
https://mathoverflow.net/questions/142192/on-a4nb4-c4nd4-and-chebyshev-polynomials

R Core Team (2013). R: A language and environment for statistical computing. R Foundation for Statistical Computing, Vienna, Austria

Roediger, P. A. (1972), "Notes on the Diophantine Equation $A^4 + aB^4 = C^4 + aD^4$", *Unpublished: Presented at the Eighteenth Conference of Army Mathematicians (Abstract and Presentation Slides available upon request)*

Roediger, P. A. (2015), "Notes on the Diophantine Equation $A^4 + aB^4 = C^4 + aD^4$", *arXiv:1510.03422*

Tomita, S. (2017), http://www.maroon.dti.ne.jp/fermat/dioph121e.html

Zajta, A. J. (1983), "Solutions of the Diophantine Equation $A^4 + B^4 = C^4 + D^4$", *Mathematics of Computation*, 41, 635-659